\documentclass{article}

\usepackage{amsmath,amssymb,amsthm}
\usepackage{url}

\def\Dbar{\leavevmode\lower.6ex\hbox to 0pt{\hskip-.23ex \accent"16\hss}D}
\def\dbar{\leavevmode\hbox to 0pt{\hskip.2ex
    \accent"16\hss}d}

\newcommand{\dd}{\mathrm{d}}
\newcommand{\T}{^\mathrm{T}}
\newcommand{\xupt}{t,\mathbf{x}(t),\mathbf{u}(t),\psi_0,\boldsymbol{\psi}(t)}
\newcommand{\xup}{t,\mathbf{x},\mathbf{u},\psi_0,\boldsymbol{\psi}}
\newcommand{\xut}{t,\mathbf{x}(t),\mathbf{u}(t)}
\newcommand{\xu}{t,\mathbf{x},\mathbf{u}}
\newcommand{\maple}{\texttt{Maple} }
\newcommand{\mapleSE}{\texttt{Maple}}

\newtheorem{theorem}{Theorem}
\newtheorem{example}{Example}
\newtheorem{definition}[theorem]{Definition}

\begin{document}

\title{Computing ODE Symmetries as Abnormal Variational Symmetries}

\author{Paulo D. F. Gouveia\\
        \texttt{pgouveia@ipb.pt}\\[0.3cm]
        Bragan\c{c}a Polytechnic Institute\\
        5301-854 Bragan\c{c}a, Portugal
        \and
        Delfim F. M. Torres\\
        \texttt{delfim@ua.pt}\\[0.3cm]
        University of Aveiro\\
        3810-193 Aveiro, Portugal
        }

\date{}

\maketitle


\begin{abstract}
We give a new computational method to obtain symmetries of
ordinary differential equations. The proposed approach appears as
an extension of a recent algorithm to compute variational
symmetries of optimal control problems [Comput. Methods Appl.
Math. {\bf 5} (2005), no.~4, pp.~387-409], and is based on the
resolution of a first order linear PDE that arises
as a necessary and sufficient condition of invariance for
abnormal optimal control problems. A computer algebra procedure
is developed, which permits to obtain ODE symmetries by the proposed method.
Examples are given, and results compared with those
obtained by previous available methods.
\end{abstract}

\smallskip

\textbf{Mathematics Subject Classification 2000:} 34-04; 49-04;
34C14; 49K15.

\smallskip


\smallskip

\textbf{Keywords.} Symmetries, variational symmetries, dynamic symmetries, ODEs, computer algebra systems, optimal control, abnormality.

\medskip


\section{Introduction}

Sophus Lie was the first to
introduce the use of symmetries into the study of differential
equations, Emmy Noether the first to recognize the important role
of symmetries in the calculus of variations.
Currently, all the computer algebra systems that address differential equations provide several tools to help the user with the
analysis of Lie symmetries. Recently, the authors developed
a computer algebra package for the automatic computation of
Noether variational symmetries in the calculus of variations
and optimal control \cite{GouveiaTorresCMAM},
now available as part of the Maple Application Center at
\url{http://www.maplesoft.com/applications/app_center_view.aspx?AID=1983}.

The omnipresent tools for Lie symmetries
provide a great help for the
search of solutions of ODEs,
their classification, order reduction,
proof of integrability, or in the construction of first integrals.
From the mathematical point of view, a ODE symmetry is described by a
group of transformations that keep the ordinary differential
equation invariant. Depending on the type of transformations one
is considering, different symmetries are obtained. An important
class of symmetries is obtained considering a one-parameter
family of transformations, which form a local Lie group. Those
transformations are often represented by a set of functions known
as the infinitesimal generators. From the practical point of view,
the determination of the infinitesimal generators that define a
symmetry for a given ODE is, in general, a complex task
\cite{Her94,Samokhin}. To address the problem, we follow a different approach.

We propose a new method for computing
symmetries of ODEs by using a  Noetherian perspective.
Making use of our previous algorithm
\cite{GouveiaTorresCMAM}, that has shown up good results
for the computation of Noether variational symmetries of problems
of the calculus of variations and optimal control,
we look to an ODE as being the control system of an optimal control problem.
Then, we obtain symmetries for the ODE
by computing the abnormal variational symmetries
of the associated optimal control problem.

This paper is organized as follows.
In \S\ref{sec:SymmOC}, the necessary concepts associated
with variational symmetries in optimal control are reviewed.
The new method for computing symmetries of ODEs
is explained in \S\ref{sec:symode}. The
method is illustrated in \S\ref{sec:exemplos}, where
we compute symmetries for three distinct ODEs
and compare the results with the ones
obtained by the standard procedures available in \mapleSE .
We end the paper with
some conclusions and comments \S\ref{sec:concl}.
The definitions of the new \maple procedure that implements our method are given in Appendix.


\section{Symmetries in optimal control}
\label{sec:SymmOC}
Without loss of generality, we consider the optimal control
problem in Lagrange form: to minimize an integral functional
\begin{equation}
\label{eq:funcionalCO} I[\mathbf{x}(\cdot),\mathbf{u}(\cdot)] =
\int_{a}^{b} L(\xut) \,\dd t
\end{equation}
subject to a control system described by a system of ordinary
differential equations of the form
\begin{equation}
\label{eq:sistCont}
\dot{\mathbf{x}}(t)= \boldsymbol{\varphi}(\xut) \, ,
\end{equation}
together with appropriate boundary conditions.
The Lagrangian $L: \mathbb{R}
\times \mathbb{R}^n \times \mathbb{R}^m \rightarrow \mathbb{R}$
and the velocity vector $\boldsymbol{\varphi}\!:\mathbb{R}
\times\mathbb{R}^n\times\mathbb{R}^m \rightarrow\mathbb{R}^n$ are
assumed to be continuously differentiable functions with respect
to all their arguments. The controls $\mathbf{u}\!:[a,
b]\rightarrow \mathbb{R}^m$ are piecewise continuous
functions; the state
variables $\mathbf{x}\!:[a, b] \rightarrow\mathbb{R}^n$
continuously differentiable functions.

The celebrated Pontryagin Maximum Principle \cite{Pontryagin62}
(PMP for short) gives a first-order necessary optimality condition. The PMP can be proved from a
general Lagrange multiplier theorem.
One introduces the Hamiltonian function
\begin{equation}
\label{eq:hamilt}
H(\xup)= \psi_0 L(\xu)+ \boldsymbol{\psi}\T   \cdot
\boldsymbol{\varphi}(\xu) \, ,
\end{equation}
where $\left(\psi_0,\boldsymbol{\psi}(\cdot)\right)$ are
the ``Lagrange multipliers'', with $\psi_0 \leq
0$ a constant and $\boldsymbol{\psi}(\cdot)$ a $n$-vectorial
piecewise $C^1$-smooth function, and the multiplier theorem
asserts that the optimal control problem is equivalent to the
maximization of the augmented functional
\begin{equation}
\label{eq:funcionalCOaugm}
J[\mathbf{x}(\cdot),\mathbf{u}(\cdot),\psi_0,\boldsymbol{\psi}(\cdot)]
=
\int_{a}^{b} \left( H(\xupt) -\boldsymbol{\psi}(t)\T
\cdot  \dot{\mathbf{x}}(t)\right)\,\dd t \, .
\end{equation}

\begin{definition}
\label{def:extr} A quadruple
$\left(\mathbf{x}(\cdot),\mathbf{u}(\cdot),\psi_0,
\boldsymbol{\psi}(\cdot)\right)$ satisfying
the Pontryagin Maximum Principle
is said to be a (Pontryagin) \emph{extremal}.
An extremal is said to be \emph{normal} when $\psi_0\neq0$,
\emph{abnormal} when $\psi_0=0$.
\end{definition}

Let $\mathbf{h}^s:[a, b] \times \mathbb{R}^n \times \mathbb{R}^m
\times \mathbb{R} \times \mathbb{R}^n \rightarrow \mathbb{R}
\times \mathbb{R}^n \times \mathbb{R}^m \times \mathbb{R}^n$ be a
one-parameter group of $\mathbb{C}^1$ transformations of the form
\begin{eqnarray}
\label{eq:transf}
\mathbf{h}^s(\xup)= \hspace{8.1cm} \nonumber\\
(h_t^s(\xup),\mathbf{h}_\mathbf{x}^s(\xup),
\mathbf{h}_\mathbf{u}^s(\xup),
\mathbf{h}_{\boldsymbol{\psi}}^s(\xup)) \, .
\end{eqnarray}
Without loss of generality, we assume that the identity
transformation of the group (\ref{eq:transf}) is obtained when the
parameter $s$ is zero:
\begin{eqnarray*}
h_t^0(\xup)=t,\;
\mathbf{h}_\mathbf{x}^0(\xup)=
\mathbf{x},\;\nonumber\\
\mathbf{h}_\mathbf{u}^0(\xup)=
\mathbf{u},\;
\mathbf{h}_{\boldsymbol{\psi}}^0(\xup)=\boldsymbol{\psi}.
\end{eqnarray*}
Associated with a one-parameter group of transformations
(\ref{eq:transf}), we introduce its \emph{infinitesimal generators}:
\begin{eqnarray}
\label{eq:transf2}
T(\xup)  &=&  \left. \frac{\partial}
{\partial{s}} h^s_t\right|_{s=0}\textrm{, }
\mathbf{X}(\xup)  = \left.
\frac{\partial}{\partial{s}} \mathbf{h}_\mathbf{x}^s\right|_{s=0}
\textrm{,}\nonumber\\
\mathbf{U}(\xup)  &=&
\left. \frac{\partial} {\partial{s}} \mathbf{h}_\mathbf{u}^s
\right|_{s=0} \textrm{, }
\boldsymbol{\Psi}(\xup)  =  \left.
\frac{\partial}{\partial{s}}
\mathbf{h}_{\boldsymbol{\psi}}^s\right|_{s=0} \textrm{. }
\end{eqnarray}

We can define variational invariance using the augmented functional (\ref{eq:funcionalCOaugm}) and the one-parameter group of transformations (\ref{eq:transf}) or
an equivalent condition in terms of the generators
\eqref{eq:transf2}:

\begin{definition}[\cite{Djukic73,Torres05}]
\label{thm:condInv} An optimal control problem
(\ref{eq:funcionalCO})-(\ref{eq:sistCont}) is said to be \emph{invariant} under \eqref{eq:transf2} or, equivalently,
\eqref{eq:transf2} is said to be a \emph{symmetry}
of the problem (\ref{eq:funcionalCO})-(\ref{eq:sistCont}) if
\begin{equation}
\label{eq:detGeradVar}
\frac{\partial H}{\partial t}T +\frac{\partial
H}{\partial \mathbf{x}}\cdot \mathbf{X} +\frac{\partial
H}{\partial \mathbf{u}}\cdot \mathbf{U} +\frac{\partial
H}{\partial \boldsymbol{\psi}}\cdot \boldsymbol{\Psi}
-\boldsymbol{\Psi}\T \cdot \dot{\mathbf{x}}
-\boldsymbol{\psi}\T \cdot \frac{\dd
\mathbf{X}}{\dd t} +H \frac{\dd T}{\dd t}=0\, ,
\end{equation}
with $H$ the Hamiltonian \eqref{eq:hamilt}.
\end{definition}

A computational algorithm to obtain the infinitesimal generators
$T$, $\mathbf{X}$, $\mathbf{U}$, and $\boldsymbol{\Psi}$
that form a variational symmetry (\ref{eq:detGeradVar})
for a given optimal control problem
(\ref{eq:funcionalCO})-(\ref{eq:sistCont}) was developed in
\cite{GouveiaTorresCMAM}. Here we remark that
the abnormal variational symmetries
(\textrm{i.e.} the ones associated with $\psi_0 = 0$)
obtained by the method introduced in \cite{GouveiaTorresCMAM}
provide symmetries for ordinary differential equations.


\section{Computing ODE symmetries from an optimal control perspective}
\label{sec:symode}

To ODE symmetries (see \textrm{e.g.} \cite{Terrab98,Zwillinger98})
we associate a smaller group of infinitesimal transformations
than that used in optimal control: we are only
interested in transformations of the independent
variable $t$ and transformations of the dependent variables $\mathbf{x}$. Let $\mathbf{g}^s:\mathbb{R} \times \mathbb{R}^{n}
\rightarrow \mathbb{R}
\times \mathbb{R}^n$ be a
one-parameter group of $\mathbb{C}^1$ transformations of the form
\begin{equation*}
\mathbf{g}^s(t,\mathbf{x})=
(g_t^s(t,\mathbf{x}),
\mathbf{g}_\mathbf{x}^s(t,\mathbf{x})) \, ,
\end{equation*}
with $\mathbf{g}^0(t,\mathbf{x})=
(t,\mathbf{x})$. Let us denote
the respective \emph{infinitesimal generators} by
\begin{equation}
\label{eq:odetransf2}
\xi(t,\mathbf{x})=\left. \frac{\partial}
{\partial{s}} g^s_t(t,\mathbf{x})\right|_{s=0}\textrm{, }
\quad \boldsymbol{\eta}(t,\mathbf{x})= \left.
\frac{\partial}{\partial{s}} \mathbf{g}_\mathbf{x}^s(t,\mathbf{x})\right|_{s=0} \,
.
\end{equation}
Our method begins by identifying a richer set of variational symmetries in the form (\ref{eq:transf2}), from which we then obtain (\ref{eq:odetransf2}). In this section, we explain in detail how to arrive to the set of generators (\ref{eq:odetransf2})
that keep an ODE invariant.

We are interested in determining symmetries for systems of
ODEs in the canonical form
\begin{gather}
\label{eq:ODEsCanon}
\begin{cases}
y_1^{(r_1)} &= \phi_1(t,y_1,\dot{y}_1\ldots y_1^{(r_1-1)},
\cdots, y_n,\dot{y}_n\ldots y_n^{(r_n-1)})\,,
\\
&\vdots\\
y_n^{(r_n)} &= \phi_n(t,y_1,\dot{y}_1\ldots y_1^{(r_1-1)},
\cdots, y_n,\dot{y}_n\ldots y_n^{(r_n-1)})\,,
\end{cases}
\end{gather}
where functions $\phi_k\!:\mathbb{R}
\times\mathbb{R}^{\Sigma_{i=1}^{n} r_i} \rightarrow\mathbb{R}$, $k = 1,\ldots,n$, are continuously differentiable with respect
to all their arguments. To write the
system \eqref{eq:ODEsCanon} of differential equations as a control system (\ref{eq:sistCont}),
we begin by converting it as a system of equations of first order. For that we introduce a new set of variables, represented by the vector $\mathbf{x}$:
\begin{eqnarray}
\label{eq:novasvar}
\mathbf{x}&=&\left[ x_1,\ldots,x_{r} \right]\T
\nonumber\\
&=&\left[y_1,\dot{y_1},\ldots,y_1^{(r_1-1)},
\ldots,
y_n,\dot{y_n},\ldots,y_n^{(r_n-1)}\right]\T \, ,
\end{eqnarray}
where $r=\Sigma_{i=1}^{n} r_i$.
With this notation, we get the control system
\begin{eqnarray}
\label{eq:SC}
 \begin{cases}
 \dot{x}_1&=x_{2}\, ,\\
 &\vdots\\
 \dot{x}_{r_1-1}&=x_{r_1}\, ,\\
 \dot{x}_{r_1}&=\phi_1(t,\mathbf{x})\, ,
\end{cases}
\cdots \,
\begin{cases}
 \dot{x}_{r_1+\ldots+r_{n-1}+1}&=x_{r_1+\ldots+r_{n-1}+2}\, ,\\
 &\vdots\\
 \dot{x}_{r_1+\ldots+r_{n}-1}&=x_{r_1+\ldots+r_{n}}\, ,\\
 \dot{x}_{r_1+\ldots+r_{n}}&=\phi_n(t,\mathbf{x})\, ,
\end{cases}
\end{eqnarray}
with $r$ state variables but no control variables (\textrm{i.e.}, (\ref{eq:SC}) is a particular case of (\ref{eq:sistCont}) where
$\boldsymbol{\varphi}$ does not depend on $\mathbf{u}$).

To use the formalism of optimal control and the notion of variational symmetry \cite{GouveiaTorresCMAM}, one thing is missing: the existence of an integral functional \eqref{eq:funcionalCO} to be minimized, and whose Lagrangian $L$ enters into the definition of the Hamiltonian (\ref{eq:hamilt}), thus being necessary for computing symmetries by  \eqref{eq:detGeradVar}.
However, if we restrict ourselves to the abnormal case, where the cost functional has no role, \textrm{i.e.} if we fix $\psi_0=0$, then the Hamiltonian $H$ does not depend on the
Lagrangian $L$ and we can look to our system \eqref{eq:SC}
as an optimal control problem. If we only consider the abnormal
case, the control system \eqref{eq:SC} is everything one needs 
to write \eqref{eq:detGeradVar} and find symmetries.

We are now in conditions to use our
\maple optimal control package \cite{GouveiaTorresCMAM} and its procedure \texttt{Symmetry} to obtain symmetries for
systems of ODEs (\ref{eq:ODEsCanon}). We only need to rewrite
(\ref{eq:ODEsCanon}) as in (\ref{eq:SC}) and then call function \texttt{Symmetry} of \cite{GouveiaTorresCMAM} for the abnormal case. After using this technique with several concrete examples,
and to be able to compare the obtained results with the ones from standard techniques, we concluded that a great manual effort is necessary at each particular problem in converting the initial system into the canonical form, then to (\ref{eq:SC}), and finally recovering the initial notation to compare the results with those obtained by the tools already available in the Computer Algebra System \mapleSE . To do all the process in an entirely automatic way, and also to optimize the algorithm, we define here a new \maple function
\texttt{odeSymm} (see Appendix) whose purpose is to compute symmetries for systems (\ref{eq:ODEsCanon}) of ODEs.

\subsection*{The algorithm}

We consider the abnormal Hamiltonian
\begin{equation*}
H(t,\mathbf{x},\boldsymbol{\psi})= \boldsymbol{\psi}\T   \cdot
\boldsymbol{\varphi}(t,\mathbf{x}) \, ,
\end{equation*}
where the velocity vector is given by
\begin{eqnarray}
\label{eq:vv}
\boldsymbol{\varphi}(t,\mathbf{x})=\left[ x_{2}, \ldots,
x_{r_1}, \phi_1(t,\mathbf{x}),
x_{r_1+2}, \ldots, x_{r_1+r_2},\phi_2(t,\mathbf{x}),\right.\cdots \nonumber\\
\cdots
,\left.x_{r_1+\ldots+r_{n-1}+2},\ldots, x_{r}, \phi_n(t,\mathbf{x})
\right]\T \, .
\end{eqnarray}
Condition \eqref{eq:detGeradVar} simplifies to
\begin{equation}
\label{eq:detGerad}
\boldsymbol{\psi}\T \cdot
\left(
\frac{\partial
\boldsymbol{\varphi}
}{\partial t} T +\frac{\partial
\boldsymbol{\varphi}
}{\partial \mathbf{x}}\cdot \mathbf{X}
-\frac{\dd
\mathbf{X}}{\dd t} +
\boldsymbol{\varphi}
\frac{\dd T}{\dd t}
\right)
+\boldsymbol{\Psi}\T \cdot \left(
\boldsymbol{\varphi} - \dot{\mathbf{x}} \right)
=0\, .
\end{equation}
Thus, given a system of ODEs,
we determine the infinitesimal generators
$\xi$ and $\boldsymbol{\eta}$ \eqref{eq:odetransf2}, which
define a symmetry of the given ODEs, in the following way:
\begin{enumerate}
\item First we rewrite the given system of ODEs in the form (\ref{eq:ODEsCanon});\footnote{We only deal with differential equations that can be written in the canonical form (\ref{eq:ODEsCanon}).}
\label{en:canon}
\item We represent the dependent variables $y_i$ and their derivatives $y_i^{(j)}$, $i=1,\ldots,n$, $j=1,\ldots,r_i-1$, present in functions $\phi_k$ of the canonical system (\ref{eq:ODEsCanon}), by a new set
    of variables $\mathbf{x}$, in accordance with \eqref{eq:novasvar};
\item We define the velocity vector $\boldsymbol{\varphi}$ \eqref{eq:vv};
\item We substitute $\boldsymbol{\varphi}$ and its partial derivatives into
\eqref{eq:detGerad};
\item From equation \eqref{eq:detGerad}, we determine the variational generators $T(t,\mathbf{x},\boldsymbol{\psi})$, $\mathbf{X}(t,\mathbf{x},
\boldsymbol{\psi})$ and $\boldsymbol{\Psi}(t,\mathbf{x},\boldsymbol{\psi})$;
\label{en:gerad}
\item In the results obtained, we go back to the initial notation of variables, by means of the inverse relations
of \eqref{eq:novasvar};
\label{en:mapinv}
\item From the obtained set of variational generators $T$, $\mathbf{X}$ and $\boldsymbol{\Psi}$, we extract the subset of generators $\xi$ and $\boldsymbol{\eta}$,
\begin{equation*}
\xi \equiv T \, , \quad \eta_i \equiv X_{1+\sum_{k=1}^{i-1} r_k}\, , \quad i=1,\ldots,n \, ,
\end{equation*}
which represent symmetries for the given system of ODEs.
\label{en:xieta}
\end{enumerate}

The infinitesimal generators obtained in step \ref{en:gerad} are
functions of the auxiliary variables $\mathbf{x}$.
Since in step~\ref{en:mapinv} the variables $\mathbf{x}$
resume to its initial meaning, we conclude that our method is able to give dynamic symmetries (\textrm{cf.} Example~\ref{ex:2}). Indeed, the generators may involve
derivatives of the dependent variables:
$$
\left(\xi,\boldsymbol{\eta}\right)\equiv
\left(\xi(t,\mathbf{y},\dot{\mathbf{y}},\ldots),
\boldsymbol{\eta}(t,\mathbf{y},\dot{\mathbf{y}},\ldots)\right)
$$
with
$(t,\mathbf{y},\dot{\mathbf{y}},\ldots)=
(t,y_1,\dot{y}_1\ldots y_1^{(r_1-1)},
\cdots, y_n,\dot{y}_n\ldots y_n^{(r_n-1)})$.

We now address the non-trivial part of our seven-step algorithm, which resides precisely in step~\ref{en:gerad}: the determination of the associated variational generators. We use the following strategy.
Expanding the total derivatives
\begin{eqnarray*}
\frac{\textrm{d}T}{\textrm{d}t} &=& \frac{\partial T}{\partial t}
+\frac{\partial T}{\partial \mathbf{x}}\cdot\dot{\mathbf{x}}
+\frac{\partial T}{\partial
\boldsymbol{\psi}}\cdot\dot{\boldsymbol{\psi}}\textrm{,}\quad \\
\frac{\textrm{d}\mathbf{X}}{\textrm{d}t} &=& \frac{\partial
\mathbf{X}}{\partial t}
+ \frac{\partial \mathbf{X}} {\partial\mathbf{x}}\cdot\dot{\mathbf{x}}
+ \frac{\partial \mathbf{X}}{\partial \boldsymbol{\psi}}\cdot
\dot{\boldsymbol{\psi}}\, \textrm{,}
\end{eqnarray*}
we write equation (\ref{eq:detGerad}) as a polynomial
\begin{equation}
\label{eq:poly}
A(t,\mathbf{x},\boldsymbol{\psi})
+B(t,\mathbf{x},\boldsymbol{\psi}) \cdot \dot{\mathbf{x}}
+ C(t,\mathbf{x},\boldsymbol{\psi})  \cdot \dot{\boldsymbol{\psi}} = 0
\end{equation}
in the $2r$ derivatives $\dot{\mathbf{x}}$ and
$\dot{\boldsymbol{\psi}}$:
\begin{eqnarray}
\label{eq:detGerad2}
\boldsymbol{\psi}\T \cdot \left(
\frac{\partial \boldsymbol{\varphi}}{\partial t} T
+\frac{\partial \boldsymbol{\varphi}}{\partial \mathbf{x}}\cdot \mathbf{X}
+\boldsymbol{\varphi} \frac{\partial T}{\partial t}
-\frac{\partial \mathbf{X}}{\partial t}
\right)
+\boldsymbol{\Psi}\T \cdot \boldsymbol{\varphi}&&
\nonumber\\
+ \left(
-\boldsymbol{\Psi}\T
+ \boldsymbol{\psi}\T \cdot \boldsymbol{\varphi} \cdot
 \frac{\partial T}{\partial\mathbf{x}}
-\boldsymbol{\psi}\T \cdot \frac{\partial \mathbf{X}}
{\partial \mathbf{x}}
\right)\cdot \dot{\mathbf{x}}&&
\\
+\left(
\boldsymbol{\psi}\T \cdot \boldsymbol{\varphi} \cdot
 \frac{\partial T}{\partial \boldsymbol{\psi}}
-\boldsymbol{\psi}\T \cdot
\frac{\partial \mathbf{X}}{\partial \boldsymbol{\psi}}
\right)
\cdot \dot{\boldsymbol{\psi}}&=0\, . & \nonumber
\end{eqnarray}
The terms in \eqref{eq:detGerad2}, which involve derivatives with
respect to vectors, are expanded in row-vectors or in matrices,
depending, respectively, if the function is a scalar function
or a vectorial one.
Equation (\ref{eq:detGerad2}) is a differential equation in the $2r+1$
unknown functions $T$, $X_1$, \dots, $X_{r}$ and
$\varPsi_1$, \dots, $\varPsi_{r}$. This equation must hold
for all $\dot{x}_1$, \dots, $\dot{x}_{r}$, $\dot{\psi}_1$,
\dots, $\dot{\psi}_{r}$, and therefore
the coefficients $A$, $B$, and $C$ of polynomial
\eqref{eq:poly} must vanish, \textrm{i.e},
\begin{equation}
\label{eq:detGerad3}
\left\{
\begin{array}{l}
\displaystyle \medskip
\boldsymbol{\psi}\T \cdot \left(
\frac{\partial \boldsymbol{\varphi}}{\partial t} T
+\frac{\partial \boldsymbol{\varphi}}{\partial \mathbf{x}}\cdot \mathbf{X}
+\boldsymbol{\varphi} \frac{\partial T}{\partial t}
-\frac{\partial \mathbf{X}}{\partial t}
\right)
+\boldsymbol{\Psi}\T \cdot \boldsymbol{\varphi}=0 \, ,\\
\displaystyle \medskip
-\boldsymbol{\Psi}\T
+ \boldsymbol{\psi}\T \cdot \boldsymbol{\varphi} \cdot
 \frac{\partial T}{\partial\mathbf{x}}
-\boldsymbol{\psi}\T \cdot \frac{\partial \mathbf{X}}
{\partial \mathbf{x}}=\mathbf{0} \, , \\
\displaystyle
\boldsymbol{\psi}\T \cdot \boldsymbol{\varphi} \cdot
 \frac{\partial T}{\partial \boldsymbol{\psi}}
-\boldsymbol{\psi}\T \cdot
\frac{\partial \mathbf{X}}{\partial \boldsymbol{\psi}}
=\mathbf{0} \, .
\end{array}
\right.
\end{equation}
Although a system of partial differential equations,
solving \eqref{eq:detGerad3} is possible using the \maple command \texttt{pdsolve} because
the system is of the first order, homogeneous,
and linear with respect to the unknown functions and their
derivatives.
We also remark that since system \eqref{eq:detGerad3} is homogeneous,
we always have, as trivial solution,
$\left(T,\mathbf{X},\mathbf{\Psi}\right) = \mathbf{0}$.

When dealing with ODEs with several dependent variables
and high-order derivatives, the number of
calculations to be done is big enough, and the help of the computer is more
than welcome. We use the computer algebra system \maple 10
to define a new procedure \texttt{odeSymm} that
does all the cumbersome computations for us -- all the steps
\ref{en:canon} to \ref{en:xieta} of our algorithm.

Our procedure \texttt{odeSymm} receives, as input, a system of ODEs, and returns, as output, a family of symmetries
$\left(\xi,\boldsymbol{\eta}\right)$ -- see definition of procedure \texttt{odeSymm} in Appendix.
To optimize the resolution of \eqref{eq:detGerad3},
we give the possibility to pass several optional parameters to
\texttt{odeSymm}. These optional parameters are described
in the Appendix and illustrated with concrete examples in \S\ref{sec:exemplos}. Here we just mention that, by default, we use the method of separation of variables (see \cite{Kythe97,Zwillinger98}) to solve \eqref{eq:detGerad3}.
More precisely, we follow \cite{Terrab95}: the generators are
replaced by the sum of unknown functions, one for each variable.
For example, $T(t,x_1,x_2) =
T_1(t)+T_2(x_1)+T_3(x_2)$. Through the optional parameters,
one can use the default solving process of the \maple solver
\texttt{pdsolve} or other specific methods (\textrm{cf.} Example~\ref{ex:1}).


\section{Illustrative examples}
\label{sec:exemplos}

To show the functionality and the usefulness of our new
procedure \texttt{odeSymm}, we consider three concrete problems found in the literature. All the examples were
carried out with \maple version 10 on a 1.4GHz 512MB RAM Pentium
Centrino. The running time of procedure \texttt{odeSymm} is
indicated, for each example, in seconds.


\begin{example}[Kamke's ODE 120]
\label{ex:1}
We begin with a first order ODE found in Kamke's book \cite{Kamke59}:
\small
\begin{verbatim}
> ode:= t*diff(y(t),t)-y(t)*(t*ln(t^2/y(t))+2)=0;
\end{verbatim}
\begin{equation*}
ode:=t{\frac {d}{dt}}y(t) -y(t)  \left( t\ln
 \left( {\frac {{t}^{2}}{y(t) }} \right) +2 \right) =0
\end{equation*}
\normalsize
To obtain symmetries of the equation we use our
\maple procedure \texttt{odeSymm}
with the additional parameter
\texttt{hint=noint}. This means that we will use
the default method of resolution of PDEs of the \maple solver \texttt{pdsolve}.
If the optional parameter \texttt{hint} is not used
(see Examples~\ref{ex:2} and \ref{ex:3} below),
our procedure \texttt{odeSymm} uses the method
of separation of variables. We obtain the following infinitesimal generators ($0.72$ sec):
\small
\begin{verbatim}
> gerad:= odeSymm(ode, y(t), split, hint=nohint);
\end{verbatim}
\[
gerad:=\left[\xi=-\frac{1}{2},\eta=-{\frac {y}{t}}\right],
\left[\xi=0,\eta=-{\frac {y}{e^{t}}}\right]
\]
\normalsize
One can test the validity of the obtained symmetries with the \texttt{symtest} command of the \texttt{DEtools} \maple package:
\small
\begin{verbatim}
> map(DEtools[symtest], [gerad], ode, y(t));
\end{verbatim}
\[
[0,0]
\]
\normalsize
The \texttt{symtest} confirm that the infinitesimal generators leave the given ODE invariant, \textrm{i.e.}, the generators obtained by our method give indeed a symmetry to Kamke's ODE~120.
It is interesting to remark that, without the knowledge of the computed symmetries, the ODE \maple solver
\texttt{dsolve} is not able to integrate the ODE:
\small
\begin{verbatim}
> dsolve(ode, y(t), class);
\end{verbatim}
\normalsize
However, when one gives to the \maple solver the infinitesimal generators found by our method, the ODE is correctly solved:
\small
\begin{verbatim}
> dsolve(ode, y(t), HINT=[gerad]);
\end{verbatim}
\[
y(t)={t}^{2}{e^{ \left( {\it C_1}-1 \right) {e^{-t}}}}
\]
\normalsize
It is also interesting to note that our method is able
to find one symmetry that is
different from the ones obtained using the standard methods of the literature. The \maple system provides nine different algorithms to compute symmetries of ODEs through the command
\texttt{symgen} of the \texttt{DEtools} package.
All the available schemes for determining the infinitesimal generators -- option \texttt{way=all} -- are not able to identify our pair of infinitesimals
$\left[\xi=0,\eta=-{\frac {y}{e^{t}}}\right]$:
\small
\begin{verbatim}
> DEtools[symgen](ode, y(t), way=all);
\end{verbatim}
\[
\left[{\it \xi}=1,{\it \eta}=2\,{\frac {y}{t}}\right]
, \left[{\it \xi}=0,{\it \eta}=y\ln  \left( {\frac {{t}^{2}}{y}}
 \right) \right]
\]
\normalsize
\end{example}


\begin{example}[Damped Harmonic Oscillator]
\label{ex:2}
We consider a harmonic oscillator with restoring force $-kx$,
emerged in a liquid in such a way that the motion of the mass $m$
is damped by a force proportional to its velocity.
Using Newton's second law one obtains, as the equation of motion,
the following second order differential equation
\cite[pp. 432--434]{LOG}:
\small
\begin{verbatim}
> EL:= m*diff(x(t),t,t)+a*diff(x(t),t)+k*x(t)=0;
\end{verbatim}
\[
EL:=m{\frac {d^{2}}{d{t}^{2}}}x(t) + a {\frac {d}{dt}}x(t) +k x(t) =0
\]
\normalsize
The symmetries for this equation
are easily obtained with our \maple procedure \texttt{odeSymm}
($1.21$ sec)
\small
\begin{verbatim}
> gerad:= odeSymm(EL, x(t), split);
\end{verbatim}
\begin{multline*}
gerad:= \left[\xi=0,\eta=x\right],
\left[\xi=0,\eta=-{\frac {m{\it x'}}{k}}\right],\left[\xi=1,
\eta=0\right],\\
\left[\xi=0,\eta={e^{-\,{\frac {at}{2m}}}}{e^{{\frac {t
\sqrt {{a}^{2}-4\,km}}{2m}}}}\right],
\left[\xi=0,\eta={e^{-{\frac {at}{2m}}}}{
e^{-{\frac {t\sqrt {{a}^{2}-4\,km}}{2m}}}}\right]
\end{multline*}
\normalsize
One can confirm that these infinitesimals represent valid symmetries for the differential equation:
\small
\begin{verbatim}
> map(DEtools[symtest], [gerad], EL, x(t));
\end{verbatim}
\[
[0, 0, 0, 0, 0]
\]
\normalsize
Note that the output of our \texttt{odeSymm} procedure includes a dynamical symmetry: the derivative of the dependent variable is present in the second pair of obtained infinitesimal generators.
\end{example}


\begin{example}[Kepler's problem]
\label{ex:3}
We now consider the Kepler's problem:
a problem of the calculus of variations -- see \cite[p.~217]{BRU}.
In this case, the Lagrangian depends on
two dependent variables $q_1$ and $q_2$:
\small
\[
L(t,\mathbf{q},\dot{\mathbf{q}})=\frac{m}{2}\left(\dot{q}_1^2
+\dot{q}_2^2\right)+\frac{K}{\sqrt{q_1^2+q_2^2}} \, .
\]
\normalsize
We will use the proposed method to determine
symmetries for the corresponding Euler-Lagrange differential equation. The Euler-Lagrange equation is trivially obtained
using our package of the calculus of variations
\cite[Example~5.2]{gouv04}:
\small
\begin{verbatim}
> L:= m/2*(v[1]^2+v[2]^2)+K/sqrt(q[1]^2+q[2]^2);
\end{verbatim}
\[
L\, := \,\frac{1}{2}\,m \left( {v_{{1}}}^{2}+{v_{{2}}}^{2} \right)
+{\frac {K}{\sqrt {{q_{{1}}}^{2}+{q_{{2}}}^{2}}}}
\]
\begin{verbatim}
> EL:= CLaws[CV][EulerLagrange](L, t, [q[1],q[2]], [v[1],v[2]]);
\end{verbatim}
\begin{multline*}
EL:=\left\{ -m{\frac {d^2}{d t^2 }}q_{{1}} (t)
-{\frac {Kq_{{1}} (t) }{ \left(  q_{{1}} (t)^{2}
+  q_{{2}} (t)  ^{2} \right) ^{3/2}}}=0,
\right. \\
\left. -m{\frac {d^2}{d t^2 }}q_{{2}} (t)
-{\frac {Kq_{{2}} (t) }{ \left(  q_{{1}} (t)^{2}
+  q_{{2}} (t) ^{2} \right) ^{3/2}}}=0 \right\}
\end{multline*}
\normalsize
In this case, the Euler-Lagrange equation is a system of two second order ODEs. Our \texttt{odeSymm} procedure
is able to determine symmetries for systems of differential equations as well ($13.32$ sec):
\small
\begin{verbatim}
> odeSymm(EL, [q[1](t),q[2](t)], split);
\end{verbatim}
\begin{multline*}
\left[\xi=0,\,\eta_{{1}}=-q_{{2}},\,\eta_{{2}}=q_{{1}}\right],
\,\left[\xi=\frac{3}{2} t,\,
\eta_{{1}}=q_{{1}},\,\eta_{{2}}=q_{{2}}\right],\,\left[\xi=1,
\,\eta_{{1}}=0,\,\eta_{{2}}=0\right]
\end{multline*}
\normalsize
It is worth to mention that this example can not be handled
by the algorithms available in \mapleSE . Indeed,
the \maple command \texttt{symgen} that
looks for a symmetry generator for a given ODE
is not able to address more than one dependent variable.
\end{example}


\section{Conclusions}
\label{sec:concl}

We have used the CAS \maple to define a new computational procedure that determines, in an automatic way, symmetries of ODEs. The automatic calculation of symmetries is
a subject much studied under
the theory of differential equations, with many results
and applications in many different areas. Our main novelty is the presentation of a new algorithm, alternative to existing ones, which looks to symmetries of ODEs as particular cases
of Noether-variational symmetries.
As explained in \S\ref{sec:symode},
our algorithm involves the resolution of a
first order, homogeneous, and linear PDE,
which is the abnormal case of the
necessary and sufficient condition of invariance
for problems of optimal control studied
with Noether's theorem
\cite{GouveiaTorresCMAM,Torres02}.
Interesting points of the proposed method are:
(i) it is based on a new approach to the subject
-- in particular, it is different from all
the nine alternative algorithms available in \mapleSE;
(ii) allows us to obtain dynamic symmetries 
for ODEs of any order;
(iii) allows to determine symmetries for systems of ODEs, when the analog \texttt{simgen} \maple command of the \texttt{DEtools} package can only obtain solutions for a single ODE.


\section*{Acknowledgements}

The authors acknowledge the support from
the \emph{Portuguese Foundation for Science and Technology} (FCT),
through the R\&D unit \emph{Centre for Research on Optimization and Control} (CEOC) of the University of Aveiro, cofinanced by
the European Community Fund FEDER/POCI 2010.
PG is also grateful to the support given by the
program PRODEP III/5.3/2003.



\section*{Appendix: the new \maple procedure \texttt{odeSymm}}

The procedure \emph{odeSymm}, introduced in
this paper, has been implemented for the computer
algebra system \maple (version~10).
The complete \maple definitions can be freely obtained from
\url{http://www.ipb.pt/~pgouveia/odeSymm.htm}
together with an online help database for the \maple system.

\begin{description}
\item[odeSymm] computes the infinitesimal generators which define
the symmetries of the ODE, or system of ODEs, specified in the
input. As explained in section \ref{sec:symode},
this procedure involves the resolution of a system of partial
differential equations.
We have used the \maple solver
\emph{pdsolve}, using, as preferential method,
the separation of the variables by sum.
\item
Output:
\begin{itemize}
\item[-] one or more lists of symmetry generators
for a given ODE, or system of ODEs
($[\xi=?,\eta_1=?,\eta_2=?,\ldots,\eta_n=?],\cdots$).
\end{itemize}
\item Syntax:
\begin{itemize}
\item[-] odeSymm(ode, x(t), opts)
\end{itemize}
\item Input:
\begin{itemize}
\item[ode -] ordinary differential equation, or a set or list of ODEs;
\item[x(t) -] any indeterminate function of one variable, or a  list of them,
representing the unknowns of the ODE problem;
\item[opts -] (optional) specify options for the \texttt{odeSymm} command,
where \emph{opts} is one or more of the following:
\begin{description}
\item[\texttt{allconst - }]
When this argument is given, the output presents all the constants given by
the \maple command \texttt{pdsolve}.
By default, that is, without option \texttt{allconst},
we eliminate redundant constants;
this is done by our \maple
procedure \texttt{reduzConst}, which is a technical routine,
and thus not provided here. Essentially, the procedure
transforms in one constant each sum of constants. The interested reader can find the \maple file
with its definition at
\url{http://www.ipb.pt/~pgouveia/odeSymm.htm}.
\item[\texttt{mindep} -]
When one wants to restrict to the minimum the dependencies
of the infinitesimal generators:
$\xi(t)$ and $\boldsymbol{\eta}(\mathbf{x})$. By default,
that is, in the absence of options
\texttt{mindep} and \texttt{alldep}, the following dependencies are considered:
$\xi(t)$ and $\boldsymbol{\eta}(t,\mathbf{x})$;
\item[\texttt{alldep} -]
All possible dependencies for the infinitesimal generators:
$\xi(t,\mathbf{x},\boldsymbol{\psi})$ and
$\boldsymbol{\eta}(t,\mathbf{x},\boldsymbol{\psi})$;
\item[\texttt{split} -]
When this argument is given, the procedure invoke the \texttt{split} command
to divide the resultant set of infinitesimal generators into uncoupled
subsets, by fixing the values for all the constants given by the \maple
command \texttt{pdsolve}. The procedure \texttt{split} is a technical routine.
The interested reader can find the \maple file
with its definition at \url{http://www.ipb.pt/~pgouveia/odeSymm.htm}.
\item[\texttt{showdep} -]
Shows, in the obtained solution, all the dependencies of the generators;
otherwise only the name of the generators is shown;
\item[\texttt{showt} -]
Shows, in the obtained solution, the dependence on the time variable
(independent variable);
otherwise, the time variable is omitted as a function parameter;
\item[\texttt{showgen} -]
Shows, in the obtained solution, besides the infinitesimal generators
$\xi$ and $\boldsymbol{\eta}$, the augmented set of variational generators,
$T$, $\mathbf{X}$ and $\boldsymbol{\Psi}$;
\item[\texttt{hint}=$<$\texttt{\textit{value}}$>$ -]
Indicate a method of solution of the PDE system (\ref{eq:detGerad3}),
where $<$\texttt{\textit{value}}$>$ is one of
 $\;\grave{}\,\texttt{+}\,\grave{}$, $\grave{}\,\texttt{*}\,\grave{}$,
or any other expression allowed by the command \texttt{pdsolve},
being also possible to use \texttt{hint=nohint} for the case one wants to use
the standard method of resolution of \mapleSE ;
by default, the system is solved by separating the variables by sum
(\texttt{hint= $\grave{}\,\texttt{+}\,\grave{}$}).
\end{description}
\end{itemize}
\end{description}
\small
\begin{verbatim}
odeSymm := proc(ODEs::{`=`, list(`=`), set(`=`)},
                                      depvars::{function,list(function)})
  local n, tt, xx, pp, k, vX, vPSI, syseqd, sol, lstGerad, valGerad, phi,
          vphi, lpsi, vpsi,Hi, t, Sr, x0, r, aux, mapx, sys, xieta, sol2;
  unprotect(Psi); unassign('T'); unassign('X'); unassign('Psi');
  unassign('psi');
  Hi:=subs(select(type,[args[3..-1]],`=`),hint);
  if Hi='hint' then Hi:=`+`; fi;
  n:=nops(depvars);
  if n=1 then x0:=[depvars] else x0:=depvars fi;
  t:=op(1,x0[1]);
  r:=[]:
  for aux in x0 do
    for k from 1 by 1 while evalb(subs(diff(aux,t$k)=_zzz,ODEs)<>ODEs) do
    od;
    r:=[r[], k-1]:
  od:
  Sr:=sum(r['i'],'i' =1..n);
  mapx:=[seq(x0[i]=_x[1+(sum(r['k'],'k'=1 ..i-1))], i=1..n)];
  mapx:=[mapx[],seq(seq(diff(x0[i],t$j)=_x[j+1+sum(r['k'],'k'=1..i-1)],
                                                   j=1..r[i]-1),i=1..n)];
  mapx:=[mapx[],seq(diff(x0[i],t$r[i])=_xx[i],i=1..n)];
  mapx:=[seq(mapx[nops(mapx)+1-i],i=1..nops(mapx))];
  sys:=subs(mapx,ODEs);
  if n=1 then solve(sys,{_xx[1]})
         else solve({sys[]},{seq(_xx[i],i=1..n)}) fi;
  phi:=subs(%, [seq(_xx[i],i=1..n)]);
  vphi:=Vector([seq([seq(_x[j],j=2+sum(r['k'],'k'=1..i-1)..sum(r['k'],
                                         'k'=1..i)), phi[i]][],i=1..n)]);
  x0:= [seq(_x[i], i = 1 .. Sr)];
  if Sr>1 then lpsi:=[seq(psi[i],i=1..Sr)] else lpsi:=[psi] fi:
  vpsi:=Vector[row](lpsi);
  if member('alldep',[args[3..-1]]) then
    tt:=t,op(x0),op(lpsi); xx:=tt; pp:=tt;
  elif member('mindep',[args[3..-1]]) then
    tt:=t; xx:=op(x0); pp:=op(lpsi);
  else tt:=t; xx:=t,op(x0); pp:=t,op(lpsi); fi:
  if Sr>1 then vX:=Vector([seq(X[i](xx), i=1..Sr)]);
          else vX:=Vector([X(xx)]); fi;
  if Sr>1 then vPSI:=Vector[row]([seq(PSI[i](pp), i=1..Sr)]);
          else vPSI:=Vector[row]([PSI(pp)]); fi;
  syseqd:={ vpsi.( map(diff,vphi,t)*T(tt)
                                 +Matrix([seq(map(diff,vphi,i),i=x0)]).vX
                          +vphi*diff(T(tt),t)-map(diff,vX,t) )+vPSI.vphi,
         convert(-vPSI+(vpsi.vphi)*Vector[row]([seq(diff(T(tt),i),i=x0)])
                     -vpsi.Matrix([seq(map(diff,vX,i),i=x0)]), 'list')[],
             convert((vpsi.vphi)*Vector[row]([seq(diff(T(tt),i),i=lpsi)])
        -vpsi.Matrix([seq(map(diff,vX,i),i=lpsi)]), 'list')[]} minus {0}:
  lstGerad:=[T(tt), convert(vX,'list')[], convert(vPSI,'list')[]];
  if Hi='nohint' then sol:=pdsolve(syseqd, lstGerad);
  else sol:=pdsolve(syseqd, lstGerad, HINT=Hi); fi;
  if not member('allconst',[args[3..-1]]) then sol:=reduzConst(sol); fi:
  valGerad:=subs(sol,lstGerad);
  sol:=[(lstGerad[i]=valGerad[i])$i=1..nops(lstGerad)];
  sol:=collect(expand(simplify(sol)),[t,op(x0),op(lpsi)]);
  if not member('showdep',[args[3..-1]]) then
    xieta:=[xi,seq(eta[i],i=1..n)];
    sol:=subs(map(i->i=op(0,i),lstGerad),sol);
  else xieta:=[xi(tt),seq(eta[i](xx),i=1..n)]; fi;
  if n=1 then xieta:=subs(eta[1]=eta,xieta) fi;
  sol:=subs('PSI'='Psi', sol);
  sol:=subs(map(i->rhs(i)=lhs(i),mapx),sol);
  xieta:=subs(map(i->rhs(i)=lhs(i),mapx),xieta);
  sol2:=[xieta[1]=rhs(sol[1]),
            seq(xieta[i+1]=rhs(sol[2+sum(r['k'],'k'=1 ..i-1)]), i=1..n)];
  if member('split',[args[3..-1]]) then sol2:=[split(sol2)]
  else sol2:=[sol2] fi;
  if member('showgen',[args[3..-1]]) then sol:=[sol,sol2[]];
  else sol:=sol2 fi;
  if n=1 then x0:=op(0,depvars) else x0:=map(i->op(0,i),depvars)[] fi;
  sol:=subs({map(i->i(t)=i,[x0,op(lpsi)])[]}, sol);
  if member('showt',[args[3..-1]]) then
    sol:=subs({map(i->i=i(t),[x0,op(lpsi)])[]},sol) fi;
  return sol[];
end proc:
\end{verbatim}
\normalsize


\end{document}